\documentclass [a4,10pt,oneside]{article}
\usepackage{amsmath, amssymb, amsthm, amsfonts, amsxtra, latexsym, amscd,
pb-diagram,  graphics,rotating,xy}

\theoremstyle{plain}
\newtheorem{dl}{Theorem}[section]
\newtheorem{bd}[dl]{Lemma}
\newtheorem{md}[dl]{Proposition}

\newtheorem{dn}[dl]{Definition}
\newtheorem{nx}[dl]{Remark}

\newtheorem{thm}{\bf Theorem}
 
\newtheorem{pro}[thm]{\bf Proposition} 
\newcommand{\A}{\mathcal{A}}
\newcommand{\tx}{\otimes }
\newcommand{\ts}{\oplus}
\newcommand{\Lh}{\frak{L}}
\newcommand{\Rh}{\frak{R}}

\newcommand{\R}{\mathcal{R}}
\newcommand{\C}{\mathcal{C}}

\begin{document}
\title{THE RELATION BETWEEN ANN-CATEGORIES AND RING CATEGORIES}
\author{Nguyen Tien Quang, Nguyen Thu Thuy and Che Thi Kim Phung\\
 {\it Hanoi National University of Education}}

\pagestyle{myheadings} 
\markboth{The relation between Ann-categories and ring categories}{N. T. Quang, N. T. Thuy and C. T. K. Phung}
\maketitle
\setcounter{tocdepth}{1}
\begin{abstract}
There are different categorizations of the definition of a {\it ring} such as {\it Ann-category} [6], {\it ring category} [2],... The main result of this paper is to prove that every axiom of the definition of a {\it ring category}, without the axiom $x_{0}=y_{0},$ can be deduced from the axiomatics of an {\it Ann-category}.
\end{abstract}
\section {Introduction}
Categories with monoidal structures $\oplus, \otimes$ (also called {\it categories with distributivity constraints}) were presented by  Laplaza [3]. M. Kapranov and V.Voevodsky [2] omitted requirements of the axiomatics of Laplaza which are related to the commutativity constraints of the operation $\otimes$ and presented the name {\it ring categories} to indicate these categories.

To approach in an other way, monoidal categories can be ``smoothed'' to become a {\it category with group structure,} when they are added the definition of invertible objects (see Laplaza [4], Saavedra Rivano [9]). Now, if the back ground category is a {\it groupoid} (i.e., each morphism is an isomorphism) then we have {\it monoidal category group-like} (see A. Fr\"{o}lich and C. T. C. Wall [1], or a {\it Gr-category} (see H. X. Sinh [11]). These categories can be classified by $H^{3}(\Pi, A)$. Each Gr-category $\mathcal G$ is determined by 3 invariants: The group $\Pi$ of classes of congruence objects, $\Pi-$module $A$ of automorphisms of the unit $1$, and an element $\overline{h}\in H^{3}(\Pi, A),$ where $h$ is induced by the associativity constraint of $\mathcal G.$ 

In 1987, in [6], N. T. Quang presented the definition of an {\it Ann-category}, as a categorization of the definition of rings, when a symmetric Gr-category (also called Pic-category) is equiped with a monoidal structure $\otimes$. In [8], [7], Ann-categories and {\it regular} Ann-categories, developed from the ring extension problem, have been classified by, respectively, Mac Lane ring cohomology [5] and Shukla algebraic cohomology [10].

The aim of this paper is to show clearly the relation between the definition of an {\it Ann-category} and a {\it ring category}.

For convenience, let us recall the definitions. Moreover, let us denote $AB$ or $A.B$ instead of $A\tx B.$
\section {Fundamental definitions}
\begin{dn} {\bf  The axiomatics of an Ann-category}\\
An Ann-category consists of:\\
i) A groupoid $\A$ together with two bifunctors $\ts,\tx:\A\times\A\longrightarrow \A.$\\
ii) A fixed object $0\in \A$ together with naturality constraints $a^+,c,g,d$ such that $(\A,\ts,a^+,c,(0,g,d))$ is a Pic-category.\\
iii) A fixed object $1\in\A$ together with naturality constraints $a,l,r$ such that $(\A,\tx,a,(1,l,r))$ is a monoidal $A$-category.\\
iv) Natural isomorphisms $\Lh,\Rh$
\[\Lh_{A,X,Y}:A\tx(X\ts Y)\longrightarrow (A\tx X)\ts(A\tx Y)\]
\[\Rh_{X,Y,A}:(X\ts Y)\tx A\longrightarrow(X\tx A)\ts(Y\tx A)\] 
such that the following conditions are satisfied:\\
(Ann-1) For each $A\in \A,$ the pairs $(L^A,\breve{L^A}),(R^A,\breve{R^A})$ determined by relations:\\
\[
\begin{aligned}
&L^A & = &A\tx- \;\;\;\; &R^A&=&-\tx A\\
&\breve{L^A}_{X,Y}& = &\Lh_{A, X, Y}\;\;\;\; &\breve{R^A}_{X,
Y}&=&\mathcal{R}_{X, Y, A}
\end{aligned}
\]
are $\ts$-functors which are compatible with $a^+$ and $c.$\\
(Ann-2) For all $ A,B,X,Y\in \A,$ the following diagrams:

\[ 
\begin{diagram}
\node{(AB)(X\ts Y)}\arrow{s,l}{\breve{L}^{AB}} \node{ A(B(X\ts
Y))}\arrow{w,t}{\quad a_{A, B, X\ts Y}\quad}\arrow{e,t}{\quad
id_A\tx\breve{L}^B\quad}\node{ A(BX\ts
BY)}\arrow{s,r}{\breve{L}^A}\\
\node{(AB)X\ts(AB)Y}\node[2]{A(BX)\ts A(BY)}\arrow[2]{w,t}
{\quad\quad a_{A, B, X}\ts a_{A, B, Y} \quad\quad}
\end{diagram}\tag{1.1} 
\]

\[ \begin{diagram}
\node{(X\ts Y)(BA)}\arrow{s,l}{\breve{R}^{BA}}\arrow{e,t}{\quad
a_{ X\ts Y, B, A}\quad}\node{((X\ts Y)B)A}\arrow{e,t}{\quad
\breve{R}^B\tx id_A\quad}\node{ (XB\ts
YB)A}\arrow{s,r}{\breve{R}^A}\\
\node{X(BA)\ts Y(BA)}\arrow[2]{e,t} {\quad\quad a_{X, B, A}\ts
a_{Y, B, A} \quad\quad}\node[2]{(XB)A\ts (YB)A}
\end{diagram}\tag{1.1'} \]
\[ \begin{diagram}
\node{(A(X\ts Y))B}\arrow{s,l}{\breve{L}^{A}\tx id_B} \node{
A((X\ts Y)B)}\arrow{w,t}{\quad a_{A, X\ts Y,
B}\quad}\arrow{e,t}{\quad id_A\tx\breve{R}^B\quad}\node{
A(XB\ts
YB)}\arrow{s,r}{\breve{L}^A}\\
\node{(AX\ts AY)B}\arrow{e,t}{\quad \breve{R}^B
\quad}\node{(AX)B\ts (AY)B}\node{A(XB)\ts A(YB)}\arrow{w,t}
{\quad a\ts a\quad}
\end{diagram}\tag{1.2} \]
\[\begin{diagram}
\node{(A\ts B)X\ts(A\ts
B)Y}\arrow{s,r}{\breve{R}^{X}\ts\breve{R}^Y} \node{(A\ts
B)(X\ts Y)}\arrow{w,t}{\breve{L}^{A\ts B}}\arrow{e,t}{
\breve{R}^{X\ts
Y}}\node{A(X\ts Y)\ts B(X\ts Y)}\arrow{s,l}{\breve{L}^A\ts\breve{L}^B}\\
\node{(AX\ts BX)\ts(AY\ts BY)}\arrow[2]{e,t} {\quad\quad v
\quad\quad}\node[2]{(AX\ts AY)\ts(BX\ts BY)}
\end{diagram}\tag{1.3} \]

\noindent commute, where $v=v_{U,V,Z,T}:(U\ts V)\ts(Z\ts T)\longrightarrow(U\ts Z)\ts (V\ts T)$ is the unique functor built from $a^+,c,id$ in the monoidal symmetric category $(\A,\ts).$\\
(Ann-3) For the unity object $1\in \A$ of the operation $\ts,$ the following diagrams:

\[
\begin{diagram}
\node{1(X\oplus Y)} \arrow[2]{e,t}{\breve{L}^1}
\arrow{se,b}{l_{X\oplus Y}}
\node[2]{1X\oplus 1Y} \arrow{sw,b}{l_X\oplus l_Y} \\
\node[2]{X\oplus Y}
\end{diagram}\tag{1.4}
\]

\[
\begin{diagram}
\node{(X\oplus Y)1} \arrow[2]{e,t}{\breve{R}^1}
\arrow{se,b}{r_{X\oplus Y}}
\node[2]{X1\oplus Y1} \arrow{sw,b}{r_X\oplus r_Y} \\
\node[2]{X\oplus Y}
\end{diagram}\tag{1.4'}
\]

commute.
\end{dn}

{\bf Remark.} The commutative diagrams (1.1), (1.1') and (1.2), respectively, mean that:\\
\[\begin{aligned} (a_{A, B, -})\;:\;& L^A.L^B &\longrightarrow \;& L^{AB}\\
(a_{-,A,B})\;:\;&R^{AB}&\longrightarrow \;&R^A.R^B\\
(a_{A, - ,B})\;:\;&L^A.R^B&\longrightarrow \;&R^B.L^A
\end{aligned}\]
are $\ts$-functors.\\
The diagram (1.3) shows that the family $(\breve{L}^Z_{X,Y})_Z=(\mathcal{L}_{-,X,Y})$ is an $\ts$-functor between the $\ts$-functors $Z\mapsto Z(X\ts Y)$ and $Z\mapsto ZX\ts ZY$, and the family $(\breve{R}^C_{A,B})_C=(\mathcal{R}_{A, B,-})$
is an $\ts$-functor between the functors $C\mapsto (A\ts B)C$ and $C\mapsto AC\ts BC.$\\
The diagram (1.4) (resp. (1.4')) shows that $l$ (resp. $r$) is an $\ts$-functor from $L^1$ (resp. $R^1$) to the unitivity functor of the $\ts$-category $\A$.
\begin{dn} {\bf  The axiomatics of a ring category}

A {\it ring category} is a category $\R$ equiped with two monoidal structures $\ts,\tx$ (which include corresponding associativity morphisms $a^{\ts}_{A,B,C},a^{\tx}_{A,B,C}$ and unit objects denoted 0, 1) together with natural isomorphisms
$$u_{A,B}:A\ts B\to B\ts A,\qquad\qquad v_{A,B,C}:A\tx(B\ts C)\to (A\tx B)\ts(A\tx C)$$
$$w_{A,B,C}:(A\ts B)\tx C\to (A\tx C)\ts(B\tx C),$$
$$x_A:A\tx 0\to 0,\qquad y_A:0\tx A\to 0.$$
These isomorphisms are required to satisfy the following conditions.

$K1 (\bullet\ts\bullet)$ The isomorphisms $u_{A,B}$ define on $\R$ a structure of a symmetric monoidal category, i.e., they form a braiding and $u_{A,B}u_{B,A}=1.$

$K2 (\bullet\tx(\bullet\ts\bullet))$ For any objects $A,B,C$ the diagram
{\scriptsize
\[
\begin{diagram}
\node{A\tx(B\ts C)}\arrow{s,l}{A\tx u_{B,C}} \arrow{e,t}{v_{A,B,C}}\node{(A\tx B)\ts(A\tx C)}\arrow{s,r}{u_{A\tx B,A\tx C}} \\
\node{A\tx(C\ts B)}\arrow{e,t}{v_{A,C,B}}\node{(A\tx C)\ts(A\tx B)}
\end{diagram}
\]}
is commutative.

$K3 ((\bullet\ts\bullet)\tx\bullet)$ For any objects $A,B,C$ the diagram
{\scriptsize
\[
\begin{diagram}
\node{(A\ts B)\tx C}\arrow{s,l}{u_{A,B}\tx C} \arrow{e,t}{w_{A,B,C}}\node{(A\tx C)\ts(B\tx C)}\arrow{s,r}{u_{A\tx C,B\tx C}} \\
\node{(B\ts A)\tx C}\arrow{e,t}{w_{B,A,C}}\node{(B\tx C)\ts(A\tx C)}
\end{diagram}
\]}
is commutative.

$K4 ((\bullet\ts\bullet\ts\bullet)\tx\bullet)$ For any objects $A,B,C,D$ the diagram
{\scriptsize
\[
\begin{diagram}
\node{(A\ts (B\ts C)D)}\arrow{s,l}{a^{\ts}_{A,B,C}\tx D} \arrow{e,t}{w_{A,B\ts C,D}}\node{AD\ts((B\ts C)D)}\arrow{e,t}{AD\ts w_{B,C,D}}\node{AD\ts(BD\ts CD)}\arrow{s,r}{a^{\ts}_{AD,BD,CD}} \\
\node{((A\ts B)\ts C)D}\arrow{e,t}{w_{A\ts B,C,D}}\node{(A\ts B)D\ts CD}\arrow{e,t}{w_{A,B,D}\ts CD}\node{(AD\ts BD)\ts CD}
\end{diagram}
\]}
is commutative.

$K5 (\bullet\tx(\bullet\ts\bullet\ts\bullet))$ For any objects $A,B,C,D$ the diagram
{\scriptsize
\[
\begin{diagram}
\node{A(B\ts (C\ts D))}\arrow{s,l}{A\tx a^{\ts}_{B,C,D}} \arrow{e,t}{v_{A,B,C\ts D}}\node{AB\ts A(C\ts D)}\arrow{e,t}{AB\ts v_{A,C,D}}\node{AB\ts(AC\ts AD)}\arrow{s,r}{a^{\ts}_{AB,AC,AD}} \\
\node{A((B\ts C)\ts D)}\arrow{e,t}{v_{A,B\ts C,D}}\node{A(B\ts C)\ts AD}\arrow{e,t}{v_{A,B,C}\ts AD}\node{(AB\ts AC)\ts AD}
\end{diagram}
\]}
is commutative.

$K6 (\bullet\tx\bullet\tx(\bullet\ts\bullet))$ For any objects $A,B,C,D$ the diagram
{\scriptsize
\[
\begin{diagram}
\node{A(B(C\ts D))}\arrow{s,l}{a^{\tx}_{A,B,C\ts D}} \arrow{e,t}{A\tx v_{B,C,D}}\node{A(BC\ts BD)}\arrow{e,t}{v_{A,BC,BD}}\node{A(BC)\ts A(BD)}\arrow{s,r}{a^{\tx}_{A,B,C}\ts a^{\tx}_{A,B,D}} \\
\node{(AB)(C\ts D)} \arrow[2]{e,t}{v_{AB,C,D}}\node[2]{(AB)C\ts (AB)D}
\end{diagram}
\]}
is commutative.

$K7 ((\bullet\ts\bullet)\tx\bullet\tx\bullet)$ Similar to the above.

$K8 (\bullet\tx(\bullet\ts\bullet)\tx\bullet)$ Similar to the above.

$K9 ((\bullet\ts\bullet)\tx(\bullet\ts\bullet))$ For any objects $A,B,C,D$ the diagram
\begin{center}
\setlength{\unitlength}{1cm}
\begin{picture}(12,4.5)
{\scriptsize
\put(0,0){$((AC\ts BC)\ts AD)\ts BD$}
\put(0.1,1.5){$(AC\ts BC)\ts(AD\ts BD)$}
\put(0.3,3){$(A\ts B)C\ts(A\ts B)D$}
\put(0.8,4.5){$(A\ts B)(C\ts D)$}

\put(10,0){$(AC\ts(BC\ts AD))\ts BD$}
\put(10,1.5){$(AC\ts(AD\ts BC))\ts BD$}
\put(9.8,3){$((AC\ts AD)\ts BC)\ts BD$}
\put(9.8,4.5){$(AC\ts AD)\ts(BC\ts BD)$}
\put(5.2,4.5){$A(C\ts D)\ts B(C\ts D)$}

\put(1.8,4.3){\vector(0,-1){0.8}}
\put(1.8,2.8){\vector(0,-1){0.8}}
\put(1.8,1.3){\vector(0,-1){0.8}}

\put(11.6,4.3){\vector(0,-1){0.8}}
\put(11.6,1.9){\vector(0,1){0.8}}
\put(11.6,1.3){\vector(0,-1){0.8}}

\put(3.1,4.6){\vector(1,0){1.9}}
\put(8.3,4.6){\vector(1,0){1.3}}
\put(9.8,0.1){\vector(-1,0){6.2}}
}

\end{picture}
\end{center}
is commutative (the notation for arrows have been omitted, they are obvious).

$K10 (0\tx 0)$ The maps $x_0,y_0:0\tx 0\to 0$ coincide.

$K11 (0\tx(\bullet\ts\bullet))$ For any objects $A,B$ the diagram
{\scriptsize
\[
\begin{diagram}
\node{0\tx(A\ts B)}\arrow{s,l}{y_{A\ts B}} \arrow{e,t}{v_{0,A,B}}\node{(0\tx A)\ts(0\tx B)}\arrow{s,r}{y_a\ts y_B}\\
\node{0}\node{0\ts 0}\arrow{w,t}{l^{\ts}_0=r^{\ts}_0}
\end{diagram}
\]}
is commutative.

$K12 ((\bullet\ts\bullet)\tx 0)$ Similar to the above.

$K13 (0\tx 1)$ The maps $y_1,r^{\tx}_0:0\tx 1\to 0$ coincide.

$K14 (1\tx 0)$ Similar to the above.

$K15 (0\tx\bullet\tx\bullet)$ For any objects $A,B$ the diagram
{\scriptsize
\[
\begin{diagram}
\node{0\tx(A\tx B)}\arrow{s,l}{y_{A\tx B}} \arrow{e,t}{a^{\tx}_{0,A,B}}\node{(0\tx A)\tx B}\arrow{s,r}{y_A\tx B}\\
\node{0}\node{0\tx B}\arrow{w,t}{y_B}
\end{diagram}
\]}
is commutative.

$K16 (\bullet\tx 0\tx\bullet),(\bullet\tx\bullet\tx 0)$ For any objects $A,B$ the diagrams
{\scriptsize
\[
\begin{diagram}
\node{A\tx(0\tx B)}\arrow{s,l}{A\tx y_{B}} \arrow[2]{e,t}{a^{\tx}_{A,0,B}}\node[2]{(A\tx 0)\tx B}\arrow{s,r}{x_A\tx B}\\
\node{A\tx 0}\arrow{e,t}{x_A}\node{0}
\node{0\tx B}\arrow{w,t}{y_B}
\end{diagram}
\]}
{\scriptsize
\[
\begin{diagram}
\node{A\tx (B\tx 0)}\arrow{s,l}{A\tx x_B}\arrow{e,t}{a^{\tx}_{A,B,0}}\node{(A\tx B)\tx 0}\arrow{s,r}{x_{A\tx B}}\\
\node{A\tx 0}\arrow{e,t}{x_A}\node{0}
\end{diagram}
\]}
are commutative.
 
$K17 (\bullet(0\ts \bullet))$ For any objects $A,B$ the diagram
{\scriptsize
\[
\begin{diagram}
\node{A\tx(0\ts B)}\arrow{s,l}{A\tx l^{\ts}_B} \arrow{e,t}{v_{A,0,B}}\node{(A\tx 0)\ts (A\tx B)}\arrow{s,r}{x_A\ts(A\tx B)}\\
\node{A\tx B}\node{0\ts(A\tx B)}\arrow{w,t}{l^{\ts}_{A\tx B}}
\end{diagram}
\]}
is commutative.

$K18 ((0\ts\bullet)\tx\bullet),(\bullet\tx(\bullet\ts 0)),((\bullet\ts 0)\tx\bullet)$ Similar to the above.
\end{dn}

\section{The relation between an Ann-category and a ring category}
In this section, we will prove that the axiomatics of a ring category, without K10, can be deduced from the one of an Ann-category. First, we can see that, the functor morphisms $a^{\oplus}, a^{\otimes}, u, l^{\ts}, r^{\ts}, v, w,$ in Definiton 2 are, respectively, the functor morphisms $a_{+}, a, c, g, d, \Lh, \Rh$ in Definition 1. Isomorphisms $x_A, y_A$ coincide with isomorphisms $\widehat{L}^A, \widehat{R}^A$ referred in Proposition 1.

We now prove that diagrams which commute in a ring category also do in an Ann-category.


K1 obviously follows from (ii) in the definition of an Ann-category.

The commutative diagrams $K2, K3, K4, K5$ are indeed the compatibility of functor isomorphisms $(L^A, \Breve L^A), (R^A, \Breve R^A)$ with the constraints $a_{+}, c$ (the axiom Ann-1).

The diagrams $K5-K9,$ respectively, are indeed the ones in (Ann-2). Particularly, K9 is indeed the decomposition of (1.3) where the morphism $v$ is replaced by its definition diagram:
 {\scriptsize\[
\begin{diagram}
\node{(P\ts Q)\ts(R\ts S)}\arrow{s,l}{v} \arrow{e,t}{a_{+}}\node{((P\ts Q)\ts R)\ts S}\node{(P\ts (Q\ts R))\ts S}\arrow{w,t}{a_{+}\ts S}\arrow{s,r}{(P\ts c)\ts S}\\
\node{(P\ts R)\ts(Q\ts S)}\arrow{e,t}{a_{+}}\node{((P\ts R)\ts Q)\ts S}\node{(P\ts(R\ts Q))\ts S.}\arrow{w,t}{a_{+}\ts S}
\end{diagram}
\]}

{\bf The proof for K17, K18}
 
\begin{bd} Let $P,$ $P^{'}$ be Gr-categories, $(a_{+}, (0, g, d)), (a^{'}_{+}, (0^{'}, g^{'}, d^{'}))$ be respective constraints, and $(F, \Breve F):P\rightarrow P^{'}$ be $\oplus$-functor which is compatible with $(a_{+}, a^{'}_{+}).$ Then $(F, \Breve F)$ is compatible with the unitivity constraints $(0, g, d)), (0^{'}, g^{'}, d^{'})).$
\end{bd}
First, the isomorphism $\widehat{F}:F0\to 0'$ is determined by the composition
\[\begin{diagram}
\node{u=F0\ts F0}\node{F(0\ts 0)}
\arrow {w,t}{\widetilde{F}}
\arrow{e,t}{F(g)}
\node{F0}\node{0'\ts F0.}\arrow{w,t}{g'}
\end{diagram}\]
Since $F0$ is a regular object, there exists uniquely the isomorphism $\widehat{F}:F0\to 0'$ such that $\widehat{F}\ts id_{F0}=u.$ Then, we may prove that $\widehat{F}$ satisfies the diagrams in the definition of the compatibility of the $\ts$-functor $F$ with the unitivity constraints.
\begin{pro}
In an Ann-category $\A,$ there exist uniquely isomorphisms
$$
\hat L^A: A\tx 0 \longrightarrow 0, \qquad \hat R^A:  0\tx A \longrightarrow 0
$$
such that the following diagrams
{\scriptsize
\[\begin{diagram}
\node{AX}\node{A(0\ts X)}\arrow{w,t}{L^A(g)}\arrow{s,r}{\breve L^A\qquad(2.1)}
\node{AX}\node{A(X\ts 0)}\arrow{w,t}{L^A(d)}\arrow{s,r}{\breve L^A\qquad(2.1')}\\
\node{0\ts AX}\arrow{n,l}{g}\node{A0\ts AX}\arrow{w,t}{\hat L^A\ts id}
\node{AX\ts 0}\arrow{n,l}{d}\node{AX\ts A0}\arrow{w,t}{id\ts \hat L^A}
\end{diagram}\]}
{\scriptsize\[\begin{diagram}
\node{AX}\node{(0\ts X)A}\arrow{w,t}{R^A(g)}\arrow{s,r}{\breve R^A\qquad(2.2)}
\node{AX}\node{(X\ts 0)A}\arrow{w,t}{R^A(d)}\arrow{s,r}{\breve R^A\qquad(2.2')}\\
\node{0\ts AX}\arrow{n,l}{g}\node{0A\ts XA}\arrow{w,t}{\hat R^A\ts id}
\node{AX\ts 0}\arrow{n,l}{d}\node{XA\ts 0A}\arrow{w,t}{id\ts \hat R^A}
\end{diagram}\]}
commute, i.e., $L^A$ and $R^A$ are U-functors respect to the operation $\ts$.
\end{pro}
\begin{proof}
Since $(L^A, \breve L^A)$ are $\ts$-functors which are compatible with the associativity constraint $a^{\ts}$ of the Picard category $(\A,\ts),$ it is also compatible with the unitivity constraint
$(0,g,d)$ thanks to Lemma 1. That means there exists uniquely the isomorphism $\hat L^A$ satisfying the diagrams $(2.1)$ and $(2.1')$. The proof for $\hat R^A$ is similar.
The diagrams commute in Proposition 1 are indeed K17, K18.
\end{proof}

{\bf The proof for 15, K16}

\begin{bd} Let $(F,\breve F), (G,\breve G)$ be $\ts$-functors between $\ts$-categories $\C, \C'$ which are compatible with the constraints $(0, g, d), (0', g', d')$ and  $\widetilde F: F(0)\longrightarrow 0', \widetilde G: G(0)\longrightarrow 0'$ are respective isomorphisms. If $\alpha: F \longrightarrow G$ in an $\ts$-morphism such that $\alpha_0$ is an isomorphism, then the diagram
{\scriptsize\[\begin{diagram}
\node{F0}\arrow[2]{r,t}{\alpha_0}\arrow{se,b}{\hat F}\node[2]{G0}\arrow{sw,b}{\hat G}\\
\node[2]{0'}
\end{diagram}\]}
commutes.
\end{bd}
\begin{proof}
Let us consider the diagram
\begin{center}
\setlength{\unitlength}{1cm}
\begin{picture}(9,3.8)
{\scriptsize
\put(-0.2,0.9){$F0$}
\put(2,0.9){$F(0\ts 0)$}
\put(4.8,0.9){$G(0\ts 0)$}
\put(7.8,0.9){$G0$}

\put(-0.4,2.5){$0'\ts F0$}
\put(2,2.5){$F0\ts F0$}
\put(4.8,2.5){$G0\ts G0$}
\put(7.6,2.5){$0'\ts G0$}

\put(3.8,0.1){$u_0$}
\put(3.6,3.7){$id\ts u_0$}
\put(0.8,1.1){$F(g)$}
\put(3.6,1.1){$u_{0\ts 0}$}
\put(6.5,1.1){$G(g)$}
\put(0.9,2.7){$\breve{F}\ts id$}
\put(3.5,2.7){$u_0\ts u_0$}
\put(6.3,2.7){$\breve{G}\ts id$}

\put(-0.3,1.7){$g'$}
\put(2.3,1.7){$\widetilde{F}$}
\put(5.1,1.7){$\widetilde{G}$}
\put(7.7,1.7){$g'$}

\put(0,0.8){\line(0,-1){0.8}}
\put(0,0){\line(1,0){8}}
\put(8,0){\vector(0,1){0.8}}

\put(0,3.6){\line(0,-1){0.7}}
\put(0,3.6){\line(1,0){8}}
\put(8,3.6){\vector(0,-1){0.7}}

\put(0,2.3){\vector(0,-1){1}}
\put(8,2.3){\vector(0,-1){1}}

\put(2.6,1.3){\vector(0,1){1}}
\put(5.4,1.3){\vector(0,1){1}}

\put(1.8,2.6){\vector(-1,0){1}}
\put(3.3,2.6){\vector(1,0){1.3}}
\put(6.1,2.6){\vector(1,0){1.3}}

\put(0.3,1){\vector(1,0){1.5}}
\put(3.3,1){\vector(1,0){1.3}}
\put(6.1,1){\vector(1,0){1.5}}

\put(3.8,3.2){(I)}
\put(1,1.7){(II)}
\put(3.7,1.7){(III)}
\put(6.4,1.7){(IV)}
\put(3.8,0.5){(V)}

}

\end{picture}
\end{center}
In this diagram, (II) and (IV) commute thanks to the compatibility of $\ts$-functors $(F,\breve F), (G,\breve G)$ with the unitivity constraints; (III) commutes since $u$ is a $\ts$-morphism; (V) commutes thanks to the naturality of $g'.$ Therefore, (I) commutes, i.e.,
$$\breve{G}\circ u_0\ts u_0=\breve{F}\ts u_0.$$
Since $F0$ is a regular object, $\breve{G}\circ u_0=\breve{F}.$
\end{proof}

\begin{md}
For any objects $X, Y\in \text{ob}\A$ the diagrams
\begin{align}
{\scriptsize\begin{diagram}
\node{X\tx(Y\tx 0)}\arrow{e,t}{id\tx \widehat{L}^Y}\arrow{s,l}{a}\node{X\tx 0}
\arrow{s,r}{\widehat L^X \qquad(2.3)}
\node{0\tx(X\tx Y)}\arrow{e,t}{\widehat R^{XY}}\arrow{s,l}{a}\node{0}\\
\node{(X\tx Y)\tx 0}\arrow{e,t}{\widehat L^{XY}}\node{0}
\node{(0\tx X)\tx Y}\arrow{e,t}{\widehat R^X\tx id}\node{0\tx Y}\arrow{n,r}{\widehat R^Y \qquad(2.3')}
\end{diagram}\nonumber}
\end{align}
{\scriptsize\[\begin{diagram}
\node{X\tx(0\tx Y)}\arrow[2]{e,t}{a}\arrow{s,l}{id\tx \hat R^Y}\node[2]{(X\tx 0)\tx Y}
\arrow{s,r}{\widehat L^X\tx id\qquad(2.4)}\\
\node{X\tx 0}\arrow{e,t}{\widehat L^X}\node{0}
\node{0\tx Y}\arrow{w,t}{\widehat R^Y}
\end{diagram}\]}
commute.
\end{md}
\begin{proof}
To prove the first diagram commutative, let us consider the diagram
{\scriptsize\[\begin{diagram}
\node{X\tx(Y\tx 0)}\arrow{e,t}{id\tx\hat L^Y}\arrow{s,l}{a}\arrow{se,t}{\widehat{L^X\circ L^Y}}
\node{X\tx 0}\arrow{s,r}{\hat L^X}\\
\node{(X\tx Y)\tx 0}\arrow{e,t}{\hat L^{XY}}\node{0}
\end{diagram}\]}
According to the axiom (1.1), $(a_{X, Y, Z})_Z$ is an $\ts$-morphism from the functor $L = L^X\circ L^Y$ to the functor $G = L^{XY}$. Therefore, from Lemma 2, (II) commutes. (I) commutes thanks to the determination of $\hat L$ of the composition $L = L^\circ L^Y$. So the outside commutes.

 The second diagram is proved similarly, thanks to the axiom (1.1'). To prove that the diagram (2.4) commutes, let us consider the diagram
{\scriptsize\[\begin{diagram}
\node{X\tx (0\tx Y)}\arrow[2]{e,t}{a}\arrow{s,l}{id\tx\widehat{R^Y}}\arrow{se,t}{\hat H}\node[2]{(X\tx 0)\tx Y}\arrow{s,r}{\widehat{L^X}\tx id}\arrow{sw,t}{\hat K}\\
\node{X\tx 0}\arrow{e,b}{\widehat{L^X}}\node{0}\node{0\tx Y}\arrow{w,b}{\widehat{R^Y}}
\end{diagram}\]}
where $H = L^X\circ R^Y$ and $K = R^Y\circ L^X$. Then (II) and (III) commute thanks to the determination of the isomorphisms $H$ and $K$. From the axiom (1.2), $(a_{X,Y,Z})_Z$ is an $\ts$-morphism from the functor $H$ to the functor $K$. So from Lemma 2, (I) commutes. Therefore, the outside commutes. The diagrams in Proposition 2 are indeed K15, K16.
\end{proof}
{\bf Proof for K11}

\begin{md} In an Ann-category, the diagram
{\scriptsize\[\begin{diagram}
\node{0\ts 0}\arrow{e,t}{g_0=d_0}\node{0}\\
\node{(0\tx X)\ts (0\tx Y)}\arrow{n,l}{\widehat{R}^X\ts \widehat{R}^Y}
\node{0\tx(X\ts Y)}\arrow{n,r}{\widehat{R}^{XY}\qquad(2.5)}\arrow{w,t}{\breve{L}^0}
\end{diagram}\]}
commutes.
\end{md}
\begin{proof}
Let us consider the diagram
\begin{center}
\setlength{\unitlength}{1cm}
\begin{picture}(14,7.3)
{\scriptsize
\put(1.2,0){$A(B\ts C)\ts 0(B\ts C)$}
\put(8.1,0){$A(B\ts C)\ts 0$}
\put(1,1.7){$(AB\ts AC)\ts(0B\ts 0C)$}
\put(7.9,1.7){$(AB\ts AC)\ts(0\ts 0)$}
\put(1.1,3.5){$(AB\ts 0B)\ts(AC\ts 0C)$}
\put(8.2,3.5){$(AB\ts 0)\ts(AC\ts 0)$}
\put(1.3,5.2){$(A\ts 0)B\ts(A\ts 0)C$}
\put(9,5.2){$AB\ts AC$}
\put(1.7,7){$(A\ts 0)(B\ts C)$}
\put(8.9,7){$A(B\ts C)$}

\put(2.9,0.8){$\breve{L}^A\ts \breve{L}^0$}
\put(2.9,2.9){$v$}
\put(2.9,4.5){$\breve{R}^B\ts \breve{R}^C$}
\put(2.9,6.2){$\breve{L}^{A\ts 0}$}

\put(9.8,0.8){$\breve{L}^A\ts d_0^{-1}$}
\put(9.8,2.9){$v$}
\put(9.8,4.5){$d_{AB}\ts d_{AC}$}
\put(9.8,6.2){$\breve{L}^A$}

\put(0.2,2.6){$\breve{R}^{B\ts C}$}
\put(11.7,2.6){$d$}

\put(5.5,0.3){$f'_A\ts id$}
\put(4.6,2){$(id\ts id)\ts(\widehat{R}^B\ts\widehat{R}^C)$}
\put(4.7,3.8){$(id\ts\widehat{R}^B)\ts(id\ts\widehat{R}^C)$}
\put(4.9,5.5){$(d_A\tx id)\ts(d_A\tx id)$}
\put(5.6,7.3){$d_A\tx id$}

\put(5.5,6.3){(I)}
\put(5.5,4.5){(II)}
\put(5.5,2.7){(III)}
\put(5.5,1){(IV)}
\put(0.5,4.3){(V)}
\put(11.2,4){(VI)}

\put(2.7,0.5){\vector(0,1){0.9}}
\put(2.7,3.3){\vector(0,-1){0.9}}
\put(2.7,5){\vector(0,-1){0.9}}
\put(2.7,6.8){\vector(0,-1){1.1}}

\put(9.6,0.4){\vector(0,1){1}}
\put(9.6,3.2){\vector(0,-1){0.9}}
\put(9.6,4){\vector(0,1){1}}
\put(9.6,6.8){\vector(0,-1){1.1}}

\put(4.2,0.1){\vector(1,0){3.7}}
\put(4.4,1.8){\vector(1,0){3.3}}
\put(4.4,3.6){\vector(1,0){3.6}}
\put(4.2,5.3){\vector(1,0){4.6}}
\put(3.9,7.1){\vector(1,0){4.7}}

\put(0,0.1){\vector(1,0){1}}
\put(12,7.1){\vector(-1,0){1.7}}

\put(0,0.1){\line(0,1){7}}
\put(0,7.1){\line(1,0){1.5}}
\put(12,0.1){\line(-1,0){1.9}}
\put(12,0.1){\line(0,1){7}}
}
\put(13,4){(2.6)}
\end{picture}
\end{center}
In this diagram, (V) commutes thanks to the axiom I(1.3), (I) commutes thanks to the functorial property of $\Lh;$ the outside and (II) commute thanks to the compatibility of the functors $R^{B\ts C},R^B,R^C$ with the unitivity constraint $(0,g,d);$ (III) commutes thanks to the functorial property $v;$ (VI) commutes thanks to the coherence for the ACU-functor $(L^A,\breve{L}^A).$ So (IV) commutes. Note that $A(B\ts C)$ is a regular object respect to the operation $\ts,$ so the diagram (2.5) commutes. We have K11.
\end{proof}
Similarly, we have K12.

{\bf Proof for K13, K14}
\begin{md}
In an Ann-category, we have $$\widehat{L}^1=l_0,\widehat{R}^1=r_0.$$
\end{md}
\begin{proof}
We will prove the first equation, the second one is proved similarly. Let us consider the diagram (2.7). In this diagram, the outside commutes thanks to the compatibility of $\ts$-functor $(L^1,\breve{L}^1)$ with the unitivity constraint $(0,g,d)$ respect to the operation $\ts;$ (I) commutes thanks to the functorial property of the isomorphism $l;$ (II) commutes thanks to the functorial property of $g;$ (III) obviously commutes; (IV) commutes thanks to the axiom I(1.4). So (V) commutes, i.e.,
$$\widehat{L}^1\ts id_{1.0}=l_0\ts id_{1.0}$$
Since 1.0 is a regular object respect to the operation $\ts,$ $\widehat{L}^1=l_0.$
\begin{center}
\setlength{\unitlength}{1cm}
\begin{picture}(5.5,4)
{\scriptsize
\put(0,0){$0\ts(1.0)$}
\put(4.5,0){$(1.0)\ts(1.0)$}
\put(1.5,1.3){$0\ts 0$}
\put(3.3,1.3){$0\ts 0$}
\put(3.3,2.5){$0\ts 0$}
\put(1.8,2.5){$0$}

\put(0.4,3.5){$1.0$}
\put(4.7,3.5){$1.(0\ts 0)$}

\put(4.3,0.1){\vector(-1,0){3.1}}
\put(4.5,3.6){\vector(-1,0){3.6}}
\put(0.5,0.3){\vector(0,1){3}}
\put(5.3,3.4){\vector(0,-1){3}}
\put(0.7,3.3){\vector(3,-2){1}}
\put(5,3.3){\vector(-3,-2){1}}
\put(0.7,0.3){\vector(1,1){0.9}}
\put(4.8,0.3){\vector(-1,1){0.9}}
\put(2.3,1.4){\vector(1,0){0.8}}
\put(3.6,1.6){\vector(0,1){0.8}}
\put(1.9,2.4){\vector(0,-1){0.8}}
\put(3.1,2.6){\vector(-1,0){1}}

\put(2.3,-0.3){$\widehat{L}^1\ts id$}
\put(1.3,0.6){$id\ts l_0$}
\put(3.5,0.6){$l_0\ts l_0$}
\put(2.5,0.6){$(V)$}
\put(2.5,1.1){$id$}
\put(-0.1,2){$g_{1.0}$}
\put(0.8,2){$(II)$}
\put(1.6,2){$g_0$}
\put(2.3,2){$(III)$}
\put(3.3,2){$id$}
\put(4.2,2){$(IV)$}
\put(5.4,2){$\breve{L}^1$}
\put(7.5,2){$(2.7)$}
\put(2.5,2.7){$g_0$}
\put(1.2,3.1){$l_0$}
\put(2.5,3.1){$(I)$}
\put(3.9,3.1){$l_{0\ts 0}$}
\put(1.6,3.7){$L^1(g_0)=id\tx g_0$}
}
\end{picture}
\end{center}
We have K14.

Similarly, we have K13.
\end{proof}
\begin{dn} An Ann-category $\A$ is {\it strong} if $\widehat{L}^0=\widehat{R}^0.$
\end{dn}
All the above results can be stated as follows
\begin{md} Each strong Ann-category is a ring category.
\end{md}
\begin{nx} In our opinion, in the axiomatics of a {\it ring category,} the compatibility of the distributivity constraint with the unitivity constraint $(1, l, r)$ respect to the operation $\otimes$ is necessary, i.e., the diagrams of (Ann-3) should be added.

Moreover, if the symmetric monoidal structure of the operation $\oplus$ is replaced with the symmetric categorical groupoid structure, then each ring category is an Ann-category.
\end{nx}
An open question: May the equation $\widehat{L}^0=\widehat{R}^0$ be proved to be independent in an Ann-category?
\begin{center}

\end{center}
\end{document}